\documentclass[11pt, a4paper, english, twoside]{article}
\usepackage{babel}   
\usepackage{url}     
\usepackage{amsmath}                                    
\usepackage{amssymb} 
\usepackage{amsthm}
\usepackage{fullpage}
\usepackage[small,bf,hang]{caption}        
\usepackage{subcaption}                    
\usepackage{sidecap}                       
\usepackage{hyperref}                      
\usepackage[noabbrev,nameinlink]{cleveref} 
\usepackage{autonum}                       
\usepackage{url}                           
\usepackage{cite} 
\usepackage[ruled,vlined]{algorithm2e}
\newtheorem{definition}{Definition}[section]

\newcommand{\G}{\mathtt{G}}
\begin{document}
\title{Symplectic Elimination}
\author{Ayan Mahalanobis\thanks{IISER Pune, Pune, INDIA. Email: ayan.mahalanobis@gmail.com}}
\date{}
\maketitle
\begin{abstract}
We develop the symplectic elimnation algorithm. This algorithm using simple row operations reduce a symplectic matrix to a diagonal matrix. This algorithm gives rise to a decomposition of an arbitrary matrix into a product of a symplectic matrix and a reduced matrix. This decomposition is similar to the SR decomposition studied for a long time, which is analogous to the QR decomposition.
\end{abstract}

\section{Introduction}
Symplectic matrices are important and ubiquitous. They appear in many applications. We developed some row-column operations on symplectic matrices along with the split and the twisted orthogonal groups~\cite{bhunia19,ayan1} and split unitary groups~\cite{ayan2}. In this paper, we deal exclusively with symplectic groups. We show, by using only row operations (defined later) one can reduce a symplectic matrix to a diagonal matrix. The effect of these row operations are the same as multiplying from the left by a symplectic matrix. These symplectic matrices are elementary and are symplectic transvections. Thus the algorithm \textbf{symplectic elimination} writes a symplectic matrix as a product of elementary symplectic  matrices and a diagonal matrix. 

The SR decomposition was studied by many over the last few decades, see~\cite{hh1,hh,ang} for more details and an earlier work. The motivation for the SR algorithm is simple, it is an analogue to the QR algorithm where a symplectic matrix replaces the orthogonal matrix. It was largely studied using the Householder and Givens transformations. Here we use a \textbf{very different approach}. We use the symplectic elimination algorithm to get to an algorithm similar to the SR algorithm and we call it the ST algorithm. This method that we illustrate, describes a symplectic elimination algorithm which naturally gives rise to a factorization algorithm similar to the QR algorithm. This method is general because the split and the twisted orthogonal groups have generators and row operations similar to the symplectic group (cf.\cite[Appendix A]{bhunia19}) and so does the split unitary group~\cite{ayan2}. Thus similar elimination algorithms giving rise to similar decompositions; similar to the QR decomposition, is possible for these groups.  
\section{Basics of the Symplectic Group}
In this paper, $\mathbf{F}$ is used for a field and there will be no restrictions on $\mathbf{F}$. All matrices are over $\mathbf{F}$. In this section, we recall few facts about symplectic groups. First and foremost, a symplectic group is made of matrices of even dimension $n (=2\ell)$ and shall be denoted by $\textrm{Sp}(\mathbf{F},n)$. 
These square matrices $X$ of size $n$ satisfy the relation $^\intercal{X}JX = J$, where $^\intercal{X}$ is the transpose of the matrix $X$ and $J$ is the matrix 
$
\begin{pmatrix}
0 & \mathrm{I}_\ell\\
-\mathrm{I}_\ell & 0
\end{pmatrix}
$. Here $\mathrm{I}_\ell$ is the identity matrix of size $\ell$. When there is no scope of confusion, I shall drop the $\ell$, and denote $\mathrm{I}_\ell$ by $\mathrm{I}$.

It is a well known fact that symplectic matrices form a group under matrix multiplication. Ordinarily, we will represent a symplectic matrix $\mathtt{G}\in\textrm{Sp}(\mathbf{F},n)$ as 
\[\mathtt{G} =
\begin{pmatrix}
\mathrm{A} & \mathrm{B}\\
\mathrm{C} & \mathrm{D}
\end{pmatrix}\] where $\mathrm{A},\mathrm{B},\mathrm{C},\mathrm{D}$ are matrices of size $\ell\times\ell$ over $\mathbf{F}$ and are referred to as blocks of $\mathtt{G}$.   
\subsection{Generators of a Symplectic Group $\textrm{Sp}(\mathbf{F},n)$}
The generators of a symplectic group, \textbf{often referred to as elementary matrices}, are symplectic matrices of two different kinds. One, that works like the usual Gaussian elimination on the blocks (Equation~\ref{eqn1}). The other kind is bit different as it works only on the lower or the upper set of blocks (Equation~\ref{eqn2}). These generators were used in our earlier work ~\cite{ayan1,ayan2,bhunia19};
where we used both row and column operations to reduce a symplectic group to an identity matrix. In this work, I develop an algorithm that uses \textbf{row operations to reduce a symplectic matrix to a diagonal matrix}. This algorithm also gives us a natural way to decompose a matrix into a product of a symplectic matrix and a simple matrix that we call a \textbf{reduced matrix}. This decomposition fail sometimes and we talk about it later. This decomposition is similar to the SR decomposition and we call it the ST decomposition.

Generators of the \textbf{first kind} are of the form
\begin{equation}\label{eqn1}
\mathtt{E}_{i,j}(\alpha)=
\begin{pmatrix}
\mathrm{R} & 0\\
0 & ^\intercal{\mathrm{R}}^{-1}
\end{pmatrix}
\end{equation} where $\mathrm{R}$ is of the for $\mathrm{I} + \alpha e_{i,j}$ for $1\leq i\neq j\leq \ell$. Here $e_{i,j}$ is a matrix of size $\ell\times\ell$ that has $\alpha$ at the $(i,j)$ position and $0$ everywhere else. 

The effect of multiplying $\mathtt{G}$ by this matrix from the left, i.e., pre-multiply $\mathtt{G}$, is the following:
\[
\begin{pmatrix}
\mathrm{R} & 0\\
0 & ^\intercal{\mathrm{R}}^{-1}
\end{pmatrix}
\begin{pmatrix}
\mathrm{A} & \mathrm{B}\\
\mathrm{C} & \mathrm{D}
\end{pmatrix} = 
\begin{pmatrix}
\mathrm{RA} & \mathrm{RB}\\
^\intercal{\mathrm{R}}^{-1}\mathrm{C} & ^\intercal{\mathrm{R}}^{-1}\mathrm{D}
\end{pmatrix}\] which is tantamount of doing two row operations on $\mathtt{G}$ simultaneously. One is acting on the upper blocks and the other on the lower blocks. The operations are reversed, i.e., if it was the $i\textsuperscript{th}$ row multiplied by $\alpha$ added to the $j\textsuperscript{th}$ row in the top blocks; in the lower blocks, it is the $j\textsuperscript{th}$ row multiplied by $-\alpha$ and added to the $i\textsuperscript{th}$ row. One can replicate this argument when $\mathtt{G}$ is multiplied by $\mathtt{E}_{i,j}(\alpha)$ on the right. In that case, a row operation will be replaced by a column operation.

Generators of the \textbf{second kind} are of the form:
\begin{equation}\label{eqn2}
\mathtt{F}_{i,j}(\alpha) =
\begin{pmatrix}
\mathrm{I} & \mathrm{R}\\
0 & \mathrm{I}
\end{pmatrix} \; \textrm{and} \;
\mathtt{H}_{i,j}(\alpha) = 
\begin{pmatrix}
\mathrm{I} & 0\\
\mathrm{R} & \mathrm{I}
\end{pmatrix}
\end{equation}
where $\mathrm{R}$ is of the form $\alpha(e_{i,j}+e_{j,i})$ where $1\leq i < j\leq \ell$ or $\alpha e_{i,i}$. Multiplying $\mathtt{G}$ from the left by $\mathtt{F}_{i,j}(\alpha)$  gives us the following equation:
\begin{equation}\label{eqn_H}
\begin{pmatrix}
\mathrm{I} & \mathrm{R}\\
0 & \mathrm{I}
\end{pmatrix}
\begin{pmatrix}
\mathrm{A} & \mathrm{B}\\
\mathrm{C} & \mathrm{D}
\end{pmatrix}
 = 
\begin{pmatrix}
\mathrm{A+RC} & \mathrm{B+RD}\\
\mathrm{C} & \mathrm{D}
\end{pmatrix}\end{equation}
and multiplying $\mathtt{G}$ from the left by $\mathtt{H}_{i,j}(\alpha)$ gives us the following equation:
\begin{equation}\label{eqn_F}
\begin{pmatrix}
\mathrm{I} & 0\\
\mathrm{R} & \mathrm{I}
\end{pmatrix}
\begin{pmatrix}
\mathrm{A} & \mathrm{B}\\
\mathrm{C} & \mathrm{D}
\end{pmatrix}
 = 
\begin{pmatrix}
\mathrm{A} & \mathrm{B}\\
\mathrm{RA+C} & \mathrm{RB+D}
\end{pmatrix}.\end{equation}
Looking at $\mathrm{R} = \alpha(e_{i,j}+e_{j,i})$, we see that in $\mathrm{RA}$ all rows except the $i\textsuperscript{th}$ and the $j\textsuperscript{th}$ are zero and the new $i\textsuperscript{th}$ row is $\alpha$ times the $j\textsuperscript{th}$ row and the new $j\textsuperscript{th}$ row is $\alpha$ times the $i\textsuperscript{th}$ row. There are two row operations that happens simultaneously; except $i=j$, when there is only one row operation of scalar multiplication of the row. This whole argument can be replicated when we multiply $\mathtt{G}_{i,j}$ from the right by $\mathtt{F}_{i,j}$ or $\mathtt{H}_{i,j}$. Then a row operation is replaced by a column operation.

In the next section, I describe the symplectic elimination algorithm; where we will need to do a row-exchange. This is the row-exchange of $i\textsuperscript{th}$ row with the $(i+l)\textsuperscript{th}$ row. This can be achieved by pre-multiplying
$\mathtt{G}$ by product of three matrices as follows:
\begin{equation}
\begin{pmatrix}\label{row-ex}
\mathrm{I} & \mathrm{R}\\
0 & \mathrm{I}
\end{pmatrix}
\begin{pmatrix}
\mathrm{I} & 0\\
\mathrm{R}^\prime & \mathrm{I}
\end{pmatrix}
\begin{pmatrix}
\mathrm{I} & \mathrm{R}\\
0 & \mathrm{I}
\end{pmatrix}
\end{equation}
where $\mathrm{R}$ is the matrix $e_{i,i}$ and $\mathrm{R}^\prime = -e_{i,i}$. A note of caution, there are two different indices in this paper, one is the index of an element in $\mathtt{G}$ and the other is the index of the same element in one of the blocks $\mathrm{A}, \mathrm{B}, \mathrm{C}$ and $\mathrm{D}$.

\section{Symplectic Elimination}
I now describe an algorithm (Algorithm~\ref{sym_elimination}) that reduces a symplectic matrix $\mathtt{G}$ to a diagonal matrix by row operations. In other words, we will take a symplecic matrix $\mathtt{G}$ and multiply it \textbf{from the left} by matrices of the form 
$\mathtt{E}_{i,j}(\alpha), \mathtt{F}_{i,j}(\alpha)$ and $\mathtt{H}_{i,j}(\alpha)$ and the output will be a diagonal matrix. Needless to say, the diagonal matrix is also a symplectic matrix and there is no restriction on the field $\mathbf{F}$. 

The algorithm was motivated by the SR algorithm~\cite{hh}. The algorithm reduces columns of the input matrix one after another. The algorithm runs in cycles, starting with the first column and ending with the $\ell\textsuperscript{th}$ column. The order of reduction in the $j\textsuperscript{th}$ cycle is reducing the $j\textsuperscript{th}$ column in $\mathrm{C}$ followed by a reduction of the same column in $\mathrm{A}$, then reducing the $j\textsuperscript{th}$ column in $\mathrm{D}$ and at last reducing the same column in $\mathrm{B}$. Note that, the $j\textsuperscript{th}$ column in $\mathrm{D}$ or $\mathrm{B}$ is the $(\ell+j)\textsuperscript{th}$ column of $\mathtt{G}$. In the algorithm $\mathtt{G}_{i,j}$ refers to the element in the $(i,j)$ position of $\mathtt{G}$.

\begin{algorithm}[H]
\DontPrintSemicolon
\KwIn{A symplectic matrix $\mathtt{G}$}
\KwOut{A diagonal matrix}
\For{column $j=1$ to $\ell$}{
\textbf{Step 1:} First the $j$ column of $\mathrm{C}$ is made zero. This is done by pre-multiplying $\mathtt{G}$ by $\mathtt{H}_{i,i}(\alpha)$. Here $i$ goes from $j$ to $\ell$. Thus we are actually working on the lower triangular part of $\mathrm{C}$ including the diagonal.\;
\textbf{Step 2:} Next the $j$ column of $A$ is reduced so that the $j$ column of $\mathtt{G}$ becomes a standard column. This is done by pre-multiplying $\mathtt{G}$ by $\mathtt{E}_{i,j}(\alpha)$ where $i$ goes from $j+1$ to $\ell$. Thus we are reducing the lower triangular part of $\mathrm{A}$ using its own diagonal element. \;
\tcp*[r]{The above process, using the symmetry of a symplectic matrix automatically makes the $l+j$ row a standard row.}
\textbf{Step 3:} Next reduce the $j$ column of $\mathrm{D}$. This is done by pre-multiplying $\mathtt{G}$ by $\mathtt{E}_{i,j}(\alpha)$ where $i$ goes from $j+1$ to $\ell$. Thus we reduce the lower triangular part of $\mathrm{D}$ to zero using its own diagonal elements.  \; 
\textbf{Step 4:} Then reduce the $j$ column of $\mathrm{B}$ to zero using $\mathtt{F}_{i,j}(\alpha)$ where $i$ ranges from $j$ to $\ell$. This makes the $l+j$ column of $\mathtt{G}$ a standard column.\;
\tcp*[r]{The above process, using the symmetry of a symplectic matrix automatically makes the $j$ row a standard row.}
}
\caption{A simple description of the symplectic elimination algorithm}
\label{simple_algorithm}
\end{algorithm}
All these reductions was happening in place, i.e., the matrix $\mathtt{G}$ was being revised with each reduction. At the end, $\mathtt{G}$ is a diagonal matrix and
the output is that diagonal matrix. This is because all that survives reduction are the diagonal elements of $\mathrm{A}$ and $\mathrm{C}$, rest of the elements are either reduced to zero or became zero due to the symmetry of the symplectic matrix.

\begin{algorithm}[h!]
\DontPrintSemicolon
 \KwIn{A symplectic matrix $\mathtt{G}$}
 \KwOut{A diagonal matrix}
 \* The matirix $\mathtt{G}$ will be reduced by left multiplication \;
 \* The matrix $\mathtt{G}$ is made up of four block $\mathrm{A}$, $\mathrm{B}$, $\mathrm{C}$ and $\mathrm{D}$ as described earlier\;
 \For{$j=1$ \textbf{to} $l$}
	{
	\* Start work on C\;
	\For{$i=j$ to $l$}
	{
 		\uIf {$\mathtt{G}_{i,j} = 0$}
 		{
                $\mathtt{T}\gets$ row-switch $(i)$ 
                \tcp*[r]{See Equation~\ref{row-ex}}
                $\mathtt{G}\gets \mathtt{T}\mathtt{G}$\;
                }
            \Else{
                $\alpha \gets - \mathtt{G}_{l+i,j}/\mathtt{G}_{i,j}$
                \tcp*[r]{See Equation~\ref{eqn_F}}
                $\mathtt{T}\gets \mathtt{H}_{i,i}(\alpha)$\;
                $\mathtt{G}\gets \mathtt{T}\mathtt{G}$\;}}
\* End work on C \;               
\* Start work on A\;
	\For{$i = j+1$ \textbf{to} $l$}{
    		\If {$\mathtt{G}_{j,j} = 0$} {
                \For{$x = j+1$ \textbf{to} $l$}{ 
                    \If {$\mathtt{G}_{x,j} \neq 0$}{
                        $T \gets \mathtt{E}_{j,x}(1)$\;}}
                $\mathtt{G} \gets T\mathtt{G}$\;}  
            $\alpha \gets - \mathtt{G}_{i,j}/\mathtt{G}_{j,j}$\;
            $\mathtt{T} \gets  \mathtt{E}_{i,j}(\alpha)$\;
            $\mathtt{G} \gets \mathtt{T} \mathtt{G}$\;
                }
 \* End work on $\mathrm{A}$\;
\* Start work on  $\mathrm{D}$\;
        \For {$i=j+1$ \textbf{to} $l$}{
            $\alpha \gets \mathtt{G}_{l+i,l+j}/\mathtt{G}_{l+j,l+j}$\;
            $\mathtt{T} \gets \mathtt{E}_{j,i}(\alpha)$\;
            $\mathtt{G} \gets\mathtt{T}\mathtt{G}$\;
            }            
\* End work on $\mathrm{D}$\;
\* Start B\;
        \For {$i=j$ \textbf{to}  $l$}{
            $\alpha = - \mathtt{G}_{i,l+j}/\mathtt{G}_{l+j,l+j}$
            \tcp*[r]{See Equation~\ref{eqn_H}}
            $\mathtt{T} = \mathtt{F}_{i,j}(\alpha)$\;
            $\mathtt{G} = \mathtt{T}\mathtt{G}$\;}
\* Completes B
                }
\Return{$\mathtt{G}$}\;
\caption{Symplectic Elimination (implementation ready format)}
\label{sym_elimination}
\end{algorithm}
\subsection{Justifications of the above algorithm}
To explain the algorithm and its inner working we need to define \emph{standard column} and \emph{standard row} in a matrix. 
Recall $\mathtt{G}$ as a defined before.
\begin{definition}[Standard Column]
The $i\textsuperscript{th}$ column in the matrix $\mathtt{G}$ for $1\leq i \leq 2\ell$ is a standard column if all the entries in the the column is zero except the entry in the $i\textsuperscript{th}$ position which is non-zero.
\end{definition}
\begin{definition}[Standard Row]
The $i\textsuperscript{th}$ row in the matrix $\mathtt{G}$ is a standard row for $1\leq i\leq 2\ell$ if all the entries in the the row is zero except the entry in the $i\textsuperscript{th}$ position which is non-zero.
\end{definition}
A diagonal matrix is one in which every row and column is standard. This is probably the most relevant example in this paper because the symplectic elimination that we develop reduces a symplectic matrix to a diagonal matrix by a series of row operations.

We now prove the above algorithm, every column is reduced to a standard column by mathematical induction on columns of $\G$. We prove, it is true for $j=1$ column and then assume that it is true for the first $j-1$ columns; and then prove it for the $j$ column, where $1\leq j\leq\ell$.

When $j=1$, note that we are reducing the first column of $\mathrm{C}$ into a zero column by using generators of the form $H_{i,i}(t)$, where $i=1$ to $\ell$. For each $i$ we create one generator and pre-multiply $\G$ with that generator. This works on the $\ell+i$ row of $\G$. This can also be seen as doing a row operation on $\G$, see Equation~\ref{eqn_F}.
The field element $\alpha = -\G_{l+i,j}/\G_{i,j}$. If $\G_{i,j} = 0$, then we do a row swap (Equation~\ref{row-ex}) and then $\G_{l+i,j}$ is zero and we are done with that entry. This is the part in Algorithm~\ref{sym_elimination} which refers to work on $\mathrm{C}$.

Then we start working with $\mathrm{A}$. In this case, like with the ordinary Gaussian elimination, we use the diagonal element to reduce all the elements below that to zero. Note that, now the whole first column of $\mathrm{C}$ is zero, so we must have a non-zero element in the first column of $\mathrm{A}$. So, if the diagonal element is zero we add with the first row the row with the non-zero entry in first column using $\mathtt{E}_{1,j}(1)$ where $1<j\leq\ell$. Then use the usual Gaussian elimination algorithm using the diagonal element as a pivot to reduce all the entries beneath it to $0$. 
This Gaussian elimination is pre-multiplying $\G$ with a matrix of the form $\mathtt{E}_{i,j}(\alpha)$, where $\alpha$ is a suitable field element. This is referred to in Algorithm~\ref{sym_elimination} as work on $\mathrm{A}$.

So now we have the first column of $\G$ as a standard column. An interesting thing happens; the $(l+1)\textsuperscript{th}$ row becomes a standard row. This happens due to the symmetry in a symplectic matrix. In particular, $^\intercal{\G}J\G=J$ implies
\begin{eqnarray}
{^\intercal{\mathrm{A}}\mathrm{C}} = {^\intercal{\mathrm{C}}\mathrm{A}} & {^\intercal{\mathrm{A}}\mathrm{D}}- {^\intercal{\mathrm{C}}\mathrm{B}} =\mathrm{I}\\
{^\intercal{\mathrm{D}}\mathrm{B}} = {^\intercal{\mathrm{B}}\mathrm{D}} & {^\intercal{\mathrm{B}}\mathrm{C} -{^\intercal{\mathrm{D}}\mathrm{A}}} = -\mathrm{I} 
\end{eqnarray}
It follows that $^\intercal{\mathrm{C}}\mathrm{A}$ is a symmetric matrix. From the fact, $\mathrm{C}$ has a zero first column, the first row of $^\intercal{C}\mathrm{A}$ is zero. Thus the first column of $^\intercal{\mathrm{C}}\mathrm{A}$ is zero. Note the first column of $\mathrm{A}$ is of the form $^\intercal(\ast,0,\ldots,0)$, where $\ast$ is a non-zero field element. Thus the first column of $^\intercal{\mathrm{C}}\mathrm{A}$ is zero implies that the first row in $\mathrm{C}$ is zero.

Furthermore, the first row of $\mathrm{D}$ is of the form $(\ast,0,\ldots,0)$. Note that the first row of ${^\intercal{\mathrm{A}}\mathrm{D}}$ is a non-zero scalar times the first row of $\mathrm{D}$ and it follows from ${^\intercal{\mathrm{A}}\mathrm{D}} - {^\intercal{\mathrm{C}}\mathrm{B}} =\mathrm{I}$ that it equals $(1,0,\ldots,0)$. Thus the row of $\mathrm{D}$ is of the required form and its non-zero element is the multiplicative inverse of the non-zero element in the first row of $\mathrm{A}$.

Now we work with the $(\ell+j)$ column of $\G$. We start with the first column of $\mathrm{D}$. We do the usual Gaussian elimination using the diagonal element in the first column to make the rest of the column beneath it zero. This involves pre-multiplying with $\mathtt{E}_{j,i}(\alpha)$ where $\alpha = \G_{\ell+i,\ell+j}/\G_{\ell+j,\ell+j}$. Recall we just proved that $G_{\ell+j,\ell+j}\neq 0$. This we refer to as the work on $\mathrm{D}$ in Algorithm~\ref{sym_elimination}.

Then we reduce the first column of $\mathrm{B}$ to zero. For this we create $\mathtt{F}_{i,j}(\alpha)$ and pre-multiply $\G$ with that, where $\alpha = -\G_{i,\ell+j}/\G_{\ell+j,\ell+j}$. These pre-multiplications are row operations defined by Equation~\ref{eqn_H}.

Again an interesting thing happens. The first row of $\G$ is now a standard row. This means that the first row of $\mathrm{A}$ is of the form $(\ast,0,\ldots,0)$ and the first row of $\mathrm{B}$ is a zero row.

To see this, we look at the equation ${^\intercal{\textrm{D}}\textrm{B}} = {^\intercal{\textrm{B}}\mathrm{D}}$, which says that ${^\intercal{\textrm{D}}\mathrm{B}}$ is a symmetric matrix. The first column of ${^\intercal{\textrm{D}}\mathrm{B}}$ is a zero column because the first column of $\mathrm{B}$ is zero. Thus the first row of ${^\intercal{\textrm{D}}\mathrm{B}}$ is a zero. This implies first row of $\mathrm{B}$ is zero. This what was required.

Next, note that ${^\intercal{\mathrm{D}}\mathrm{A}} - {^\intercal{\mathrm{B}}\mathrm{C}} =\mathrm{I}$. Then look at the first row of the matrix ${^\intercal{\mathrm{D}}\mathrm{A}} - {^\intercal{\mathrm{B}}\mathrm{C}}$, the first row of ${^\intercal{\mathrm{B}}\mathrm{C}}$ is zero. Thus the first row of ${^\intercal{\mathrm{D}}\mathrm{A}}$, is of the form $(\ast,0,\ldots,0)$ as it is a scalar times the first row of $\mathrm{A}$ and equals $(1,0,\ldots,0)$. Thus the first column of $\mathrm{D}$ is of the form $^\intercal(\ast,0,\ldots,0)$. Furthermore, it shows that the non-zero entry in the first column of $\mathrm{A}$ and $\mathrm{D}$ are inverses of one another.

Thus the first cycle of the algorithm is, reduce the first column and the $(\ell+1)$ column. At the end of this cycle, we have the first row and the $(\ell+1)$ row a standard row and then the first and the $(\ell+1)$ column a standard column. Now we start the second cycle by reducing the second column and the $(\ell+2)$ column. This gives rise to the first and the second column as standard, $(\ell+1)$ and $(\ell+2)$ column as standard. Similar is the case with rows. Assume that we have been able to do this reduction for the first $j-1$ columns. Now, for all practical purposes $\mathrm{A}$, $\mathrm{B}$, $\mathrm{C}$ and $\mathrm{D}$ has changed to a $(\ell-j)\times (\ell-j)$ matrix. We run the same cycle on these matrices and place them back in the larger matrix. This completes induction. 
\section{ST decomposition}
The above algorithm (Algorithm~\ref{sym_elimination}) can be used to decompose a non-singular matrix $\mathtt{M}=\mathrm{ST}$; where $\mathrm{S}$ is a symplectic matrix and $\mathrm{T}$ is a \emph{reduced matrix}. When the \textbf{input to the algorithm is changed from a symplectic matrix to a non-singular matrix}, it computes the decomposition $\mathrm{S}$ and $\mathrm{T}$.
In other words, the matrices that we multiply from the left, if accumulated, gives rise to the symplectic matrix $\mathtt{S}$ and what is left at the end of the algorithm is the reduced matrix $\mathtt{T}$. The reason behind that is simple, the entries that become zero in the symplectic case due to its symmetry does not become zero for an arbitrary matrix.
Thus one can reduce any non-singular matrix to a reduced matrix using symplectic matrices. However, the algorithm \textbf{sometime fails} to complete and then the decomposition is not possible. This happens when there is a zero element in the diagonal in $\mathrm{D}$, i.e., $\mathtt{G}_{\ell+j,\ell+j} = 0$ for some $j$ during the algorithm. This can happen because the matrix $\mathtt{M}$ is no longer symplectic. However, when this happens, if we continue without reducing this $\ell+j$ step in the algorithm, we are left with a (sub)column in the $\ell+j$ position of $\G$ and the other columns are reduced. This algorithm is similar to the SR algorithm~\cite[Table 2.5]{hh} and so we call it the ST algorithm. Similar reduction is present in the work of Bunse-Gerstner~\cite[Page 60]{ang}.
\begin{definition}[Reduced Matrix]
A matrix $\mathrm{T}$ of even order $(2\ell)$ made up of four blocks $\mathrm{A}, \mathrm{B},\mathrm{C}$ and $\mathrm{D}$ of size $\ell\times\ell$; written as a block matrix
$
\begin{pmatrix}
\mathrm{A} & \mathrm{B}\\
\mathrm{C} & \mathrm{D}
\end{pmatrix}
$,
is a reduced matrix if $\mathrm{A}$ and $\mathrm{D}$ are upper-triangular matrix and $\mathrm{B}$ and $\mathrm{C}$ are nil-triangular (upper) matrix. A nil-triangular (upper) matrix is a upper triangular matrix in which all the elements on the diagonal are zero. 
\end{definition} 
\section{Conjugation by Symplectic Matrices}
There are many applications of conjugating (also known as similarity transformations) by symplectic matrices. Williamson's theorem~\cite{williamson} is one such example. Conjugation is often used to reduce a matrix to a condensed form, such that, it is easy to compute the eigenvalues of the condensed form. That gives the eigenvalues of the original matrix. In the case of symplectic matrices, computing eigenvalues is important, see~\cite[Chapter 2]{hh}. It will be interesting to see, if the algorithm developed here aids in computing eigenvalues of symplectic matrices. 
 
Let $M_1 = S^{-1}MS$ where $M$ is an arbitrary matrix and $S$ is a symplectic matrix. Then using symplectic elimination one can write $S=S_1S_2\cdots S_kD$ where $S_i$ are elementary symplectic matrices, $D$ a diagonal matrix and $k$ a positive integer. Then computing $M_1$ is conjugating $M$ by elementary symplectic matrices consecutively and then by a diagonal matrix. Since each $S_i$ was produced by an algorithm and is not random, there is hope for a general purpose algorithm to compute the conjugation depending on $M$. 
\bibliography{paper}
\end{document}